\documentclass[a4paper]{amsart}
\input{article.sty}

\title[The CZ(\lowercase{p}) inequality on non\hyp{}compact manifolds]{The $L^p$\hyp{} Calder\'on\hyp{}Zygmund inequality on non\hyp{}compact manifolds of positive curvature}
\author{Ludovico Marini} 
\address[L. ~Marini]{Dipartimento di Matematica e Applicazioni,
Università degli Studi di Milano-Bicocca, Via R. Cozzi 55, I-20125, Milano}
\email{l.marini9@campus.unimib.it}

\author{Giona Veronelli }
\address[G.~Veronelli]{Dipartimento di Matematica ed Applicazioni,
Università degli Studi di Milano-Bicocca, Via R. Cozzi 55, I-20125, Milano}
\email{giona.veronelli@unimib.it}

\begin{document}
\begin{abstract}
We construct, for $p>n$, a concrete example of a complete non\hyp{}compact $n$\hyp{}dimensional Riemannian manifold of positive sectional curvature which does not support any $L^p$-Calder\'on\hyp{}Zygmund inequality:
\begin{equation*}
    \forall\,\varphi\in C^{\infty}_c(M),\qquad\|\Hess \varphi \|_{L^p}\le C(\|\varphi\|_{L^p}+\|\Delta\varphi\|_{L^p}).
\end{equation*}
The proof proceeds by local deformations of an initial metric which (locally) Gromov\hyp{}Hausdorff converge to an Alexandrov space.
In particular, we develop on some recent interesting ideas by G. De Philippis and J. N\'u\~nez\hyp{}Zimbron dealing with the case of compact manifolds.
As a straightforward consequence, we obtain that the $L^p$\hyp{}gradient estimates and the $L^p$\hyp{}Calder\'on\hyp{}Zygmund inequalities are generally not equivalent, thus answering an open question in literature. 
Finally, our example gives also a contribution to the study of the (non-)equivalence of different definitions of Sobolev spaces on manifolds.
 \end{abstract}

\maketitle

\section{Introduction}
\label{sec:introduction}

Let $(M,g)$ be an $n$\hyp{}dimensional complete Riemannian manifold. 
We say that $(M,g)$ supports an $L^p$\textsl{\hyp{}Calder\'on\hyp{}Zygmund inequality} for some $p\in(1,\infty)$, if there exists a constant $C>0$ such that
\begin{equation}
\label{e: CZp}
\tag{CZ(p)}
\|\Hess \varphi \|_{L^p}\le C(\|\varphi\|_{L^p}+\|\Delta\varphi\|_{L^p}), \qquad \forall\,\varphi\in C^\infty_c(M).
\end{equation}
Here, $\Hess\varphi = \nabla^2\varphi$ denotes the second order covariant derivative tensor and $\Delta$ is the (negatively defined) Laplace-Beltrami operator of $(M,g)$; both tensorial and $L^p$ norms are computed with respect to the Riemannian metric $g$, with the common abuse of notation $\|\Hess \varphi \|_{L^p}=\| \,|\Hess \varphi|\, \|_{L^p}$.

Calder\'on\hyp{}Zygmund inequalities were first established by a work of A. Calder\'on and A. Zygmund, \cite{CZ1952}, in the Euclidean space $\R^n$, where in fact one has the stronger
\begin{equation*}
\|\Hess \varphi \|_{L^p}\le C(p, n)\|\Delta\varphi\|_{L^p}, \qquad \forall\,\varphi\in C^\infty_c(\R^n);
\end{equation*}
see also \cite[Theorem 9.9]{GT1998}.
This inequality is a fundamental tool, for instance, as an \textit{a priori} estimate in the regularity theory of elliptic PDEs.
In the Riemannian setting, $CZ(p)$ is known to hold for compact manifolds $(M,g)$. Here, the constant $C$ clearly depends on the Riemannian metric. 
The case of complete non\hyp{}compact manifolds is much less understood. In this context a systematic study of Calder\'on\hyp{}Zygmund inequalities was recently initiated by B. G\"uneysu and S. Pigola in \cite{G2016,GP2015}. 
Since then, geometric analysts have shown an increasing interest towards the topic, both concerning the (non\hyp{})existence of $CZ(p)$ on a given manifold, and the interaction of Calder\'on\hyp{}Zygmund theory with other related issues \cite{GP2018,DPNZ2019,GPM2019,IRV2019,L2020,V2020}; see in particular the very nice recent survey \cite{P2020}.

When the $C^{1,\alpha}$-harmonic radius of $(M,g)$ is positive, a computation in a harmonic coordinate system together with a covering argument allows to reduce the Riemannian problem to the Euclidean setting. 
Using this strategy, the inequality $CZ(p)$ was proved to be true in the whole range $p\in (1,\infty)$ for manifolds of bounded Ricci curvature (both from above and below) and positive injectivity radius; see \cite[Theorem C]{GP2015}. 
Note that the limit cases $CZ(1)$ and $CZ(\infty)$ are disregarded as they fail to be true even in the Euclidean space, \cite{dLM1962,O1962}. 
On the other hand, manifolds which do not support $CZ(p)$ have been recently constructed in \cite{GP2015,L2020,V2020}, hence, the validity of an $L^p$\hyp{}Calder\'on\hyp{}Zygmund inequality in the Riemannian setting needs some geometric assumptions. 
It is worth mentioning that in the cited counterexamples the Ricci curvature of the manifold at hand is always unbounded from below. 

Unsurprisingly, the $L^2$ case is a peculiar one. 
Indeed, one can use the Bochner formula to obtain a much stronger result. 
\begin{theorem}[Theorem B in \cite{GP2015}]\label{th: CZ(2)}
Let $(M, g)$ be a complete Riemannian manifold, if $\Ric \geq -K^2$ then $CZ(2)$ holds on $M$ with a constant depending only on $K$.
\end{theorem}
It is worthwhile to observe that the Theorem above is sharp in the following sense which, to the best of our knowledge, has not been yet observed so far.
\begin{mytheorem}
\label{th: sharp}
For each $m\ge 2$ and $p\in (1,\infty)$, and for each increasing function $\lambda: [0,+\infty)\to \R$ such that $\lambda(t) \to +\infty$ as $t\to \infty$, there exists a complete Riemannian manifold $(M,g)$ satisfying $\min\Sect (x) \ge -\lambda(r(x))$ for $r(x)$ large enough, and which does not support an $L^p$-Calder\'on-Zygmund inequality $CZ(p)$. Here $r(x)$ is the Riemannian distance from some fixed reference point $o\in M$ and the $\min$ is over all the sectional curvatures at the point $x$. 
\end{mytheorem}
In particular, it is not possible to obtain $CZ(2)$ under negative decreasing curvature bounds (for instance $\Ric (x) \ge -Cr^\alpha(x)$ for some $\alpha>0$), as it is the case for the closely related problem of the density of smooth compactly supported functions in the Sobolev space $W^{2,p}(M)$ (see \cite{IRV2019}). Under this milder condition, however, a disturbed $CZ(p)$ holds, \cite[Section 6.2]{IRV-preprint}.

According to Theorem \ref{th: CZ(2)} and Theorem \ref{th: sharp}, the following question naturally arises.
\begin{question}[Conjectured for $\Ric\ge 0$ in \cite{G2016}, p. 177]
    \label{q: batu}
	Suppose that $(M,g)$ is geodesically complete and has lower bounded Ricci curvature. Does $CZ(p)$ hold on $(M,g)$ for all $p\in(1,\infty)$? 
\end{question}
Strong evidence for a negative answer comes from a deep and recent result by G. De Philippis and J. N\'u\~nez\hyp{}Zimbron who proved the impossibility to have a Calder\'on\hyp{}Zygmund theory on compact manifolds with constants depending only on a lower bound on the sectional curvature, at least when $p>n$. 
Namely, when $p>n$ one can find a sequence of compact, non\hyp{}negatively curved Riemannian manifolds $\{(M_j,g_j)\}_{j=1}^\infty$ for which the best constant in $CZ(p)$ is at least $j$; see \cite[Corollary 1.3]{DPNZ2019}.

The main result of this short note gives a concrete and final answer to  \Cref{q: batu}, even under the stricter assumption of positive sectional curvature.

\begin{mytheorem}
\label{t: main}
For every $n\ge 2$ and $p>n$, there exists a complete, non\hyp{}compact $n$\hyp{}dimensional Riemannian manifold $(M, g)$ with $\Sect(M) > 0$ such that $CZ(p)$ fails. 
\end{mytheorem}

With respect to the argument in \cite{DPNZ2019}, our main contribution consists in proving the existence of a fixed Riemannian manifold on which $CZ(p)$ can not hold, whatever constant $C$ one takes (as explained above, such a result is clearly impossible in the compact setting). 

To prove \cite[Corollary 1.3]{DPNZ2019}, the authors considered a sequence of smooth non\hyp{}negatively curved $n$\hyp{}dimensional compact manifolds Gromov\hyp{}Hausdorff approaching a compact $RCD(0,n)$ space $X$ with a dense set of singular points.
A bound on the constant $C$ in $CZ(p)$  along the sequence, combined with a Morrey inequality, would imply that all functions on $X$ with Laplacian in $L^{p>n}$ are $C^1$. 
On the other hand, De Philippis and N\'u\~nez\hyp{}Zimbron proved in \cite[Theorem 1.1 ]{DPNZ2019} that the gradient of a harmonic function (or more generally of any function whose Laplacian is in $L^{p>n}$) vanishes at singular points of an $RCD$ space.
By a density argument, this would imply that all harmonic functions on $X$ are constant which is impossible.   

To achieve our result, we localize this procedure.  
The key observation is the fact that the argument is indeed local and can be repeated on infinitely many singular perturbations scattered over a non\hyp{}compact manifold.
Namely, we begin with a complete non\hyp{}compact manifold $(M,g)$ with $\Sect(M) > 0$. 
In the interior of infinitely many separated sets $\{\mathfrak D_j\}_{j\in\mathbb{N}}$ of $M$ we take sequences of local perturbations $g_{j,k}$ of the original metric $g$ such that all the $g_{j,k}$ have $\Sect > 0$ and $g_{j,k}$ Gromov\hyp{}Hausdorff converges to an Alexandrov metric $d_{j,\infty}$ on $M$ of non\hyp{}negative curvature (hence $RCD(0,n)$). 
In particular, the metric $d_{j, \infty}$ is singular on a dense subset of $\mathfrak{D}_j$.
Next, we observe that a neighborhood of each $\mathfrak D_j$ can be seen as a piece of a compact space whose metric is smooth outside $\mathfrak D_j$, so that De Philippis and N\'u\~nez\hyp{}Zimbron's strategy can be applied locally to the sequence $g_{j,k}$. 
Accordingly, we find a (large enough) $k$ and a function $v_j$ compactly supported in a small neighborhood of $\mathfrak{D}_j$ such that the following estimate holds with respect to the metric $g_{j,k}$
\begin{equation*}
    \Vert \Hess v_j \Vert_{L^p} > j \left(\Vert \Delta v_j \Vert_{L^p} + \Vert v_j \Vert_{L^p}\right).
\end{equation*}
Gluing together all the local deformations of the metric, we thus obtain a smooth manifold on which no constant $C$ makes $CZ(p)$ true.

It is worthwhile to note that, to the best of our knowledge, the problem of extending to $1<p\le n$ De Philippis and N\'u\~nez\hyp{}Zimbron's result (and thus our extension) is completely open, except for the case $p=2$ alluded to above. 
Notably, it is not known whether a lower bound on the sectional curvature suffices to have the validity of $CZ(p)$ for any $p\in (1,n]$. 
The main obstruction to reproduce the strategy detailed above is the lack of a Morrey embedding when $p\leq n$. 
Accordingly, the gradient of a harmonic functions on the singular space could be non\hyp{}continuous. 
We wonder whether an $L^q$ control on the gradient of harmonic functions, for $q$ large enough, could suffice, thus permitting to lower the threshold $n$ for $p$ in Theorem \ref{t: main}.

\vspace{\baselineskip}

\Cref{t: main} confirms the strong indications carried by \cite{DPNZ2019} and sheds light on the conditions necessary to the validity of $CZ(p)$ on complete non\hyp{}compact Riemannian manifolds. 
On the other hand, our result answers as a byproduct two other related questions.
Beyond the importance that Calder\'on\hyp{}Zygmund inequalities have in themselves, their validity has consequences on other topics in the field.
For instance, $CZ(p)$ is related to a class of functional inequalities known as $L^p$\textsl{\hyp{}gradient estimates}, i.e.
\begin{equation}\label{eq: Lp-grad}
\|\nabla \varphi \|_{L^p}\le C(\|\varphi\|_{L^p}+\|\Delta\varphi\|_{L^p})
\end{equation}
for all $\varphi\in C^\infty_c(M)$.
These gradient estimates are known to hold on any complete Riemannian manifold (actually in a stronger multiplicative form) for $p\in(1,2]$, \cite{CD2003}. 
To the best of our knowledge, it is still unknown if \eqref{eq: Lp-grad} holds as well for $p>2$ without further assumption.
Nonetheless, a Riemannian manifold supports an $L^p$\hyp{}gradient estimate whenever $CZ(p)$ holds on $M$, \cite[Corollary 3.11]{GP2015}. 
This naturally leads to the following question, raised by B. Devyver.

\begin{question}[see Section 8.1 in \cite{P2020}]\label{q: baptiste}
Are $L^p$\hyp{}gradient estimates and $L^p$\hyp{}Calder\'on\hyp{}Zygmund inequalities equivalent? 
\end{question}

Since $L^p$\hyp{}gradient estimates are known to hold when the Ricci curvature is bounded from below, \cite{CTT2018}, \Cref{t: main} gives a negative answer also to \Cref{q: baptiste}.

\begin{mycorollary}
For any $n \geq 2$ and $p> n$, there exists a complete Riemannian manifold $(M,g)$ supporting the $L^p$\hyp{}gradient estimate \eqref{eq: Lp-grad} on which $CZ(p)$ does not hold. 
\end{mycorollary}

Another important feature of $CZ(p)$ is its interaction with Sobolev spaces. 
Unlike the Euclidean setting, on a Riemannian manifold there exist several non-necessarily equivalent definitions of $k$-th order $L^p$ Sobolev space; see for instance the introduction in \cite{V2020} for a brief survey. 
The role of $CZ(p)$ in the density problem of compactly supported functions in the Sobolev space is by now well understood; see \cite[Remark 2.1]{V2020}. 
Here, we consider the spaces
\[
W^{2,p}(M)=\{f\in L^p:\ \nabla f\in L^p,\ \Hess\, f\in L^p\},
\]
and
\[
H^{2,p}(M)=\{f\in L^p:\ \Delta f\in L^p\},
\]
endowed with their canonical norms.
Here, the gradient, the Hessian and the Laplacian are interpreted in the sense of distributions. 
Note that the space $H^{2,p}$ can be interpreted as the maximal self-adjoint realization of $\Delta:C^\infty_c\to C^\infty_c$ in $L^p$. 
By definition, $W^{2,p}(M)\subset H^{2,p}(M)$. 
Moreover, if $CZ(p)$ and \eqref{eq: Lp-grad} hold on $M$ one has
\begin{equation*}
\|\nabla \varphi \|_{L^p}+\|\Hess \varphi \|_{L^p}\le C(\|\varphi\|_{L^p}+\|\Delta\varphi\|_{L^p}),\quad\forall\,\varphi\in C_c^{\infty}(M).
\end{equation*}
Thanks to a density result due to O. Milatovic  (see \cite[Appendix A]{GPM2019}), the latter estimate holds for any $\varphi\in H^{2,p}$. 
Thus, $H^{2,p}=W^{2,p}$ whenever $CZ(p)$ and \eqref{eq: Lp-grad} hold on $M$, for instance, when the geometry of $M$ is bounded. 
Conversely, examples proving that $H^{2,p}\neq W^{2,p}$ on wildly unbounded geometries are known; \cite{D1981,V2020}.
In this direction, as a corollary of the proof of \Cref{t: main}, we get the following interesting observation. 

\begin{mycorollary}
\label{cor:counterexample to sobolev inclusion}
For every $n\geq 2$ and $p>n$ there exists a complete, non\hyp{}compact $n$\hyp{}dimensional Riemannian manifold with $\Sect(M) > 0$ such that $W^{2, p}(M) \subsetneq H^{2, p}(M)$. 
\end{mycorollary}

The paper is organized as follows. 
In \Cref{sec:singular space} we construct the sequence of local deformations of the initial smooth metric, each Gromov\hyp{}Hausdorff converging to an Alexandrov space of positive curvature with a locally dense cluster of singular points. 
In \Cref{sec: Poisson} we specify to our setting a proposition by S. Pigola (which collects a series of previous results due to S. Honda) about the convergence of functions defined on a Gromov\hyp{}Hausdorff converging sequence of manifolds. 
Finally, in \Cref{sec: proof} we conclude the proofs of \Cref{t: main} and \Cref{cor:counterexample to sobolev inclusion}.

\section{The singular space and its smooth approximations}
\label{sec:singular space}
It is well known from previous literature that for every $n\ge 2$ one can always construct a compact, convex set $C \subset \R^{n+1}$ whose boundary $X = \partial C$ is an Alexandrov space with $\Curv (X) \geq 0$ and a dense set of singular points. 
The first example of such spaces is due to Y. Otsu and T. Shioya in dimension 2, \cite[Example (2)]{OS1994}, although the result holds in arbitrary dimension. 
Observe that the space $X$ can be GH approximated with a sequence of smooth manifolds $X_k$ of non\hyp{}negative sectional curvature; see the proof of \cite[Theorem 1]{AKP2008}.

In the following, we would like to localize this construction inside a compact set of a complete, non\hyp{}compact manifold. 
Indeed, we prove that a smooth and strictly convex function can always be perturbed on a compact set by introducing a dense sequence of singular points. 
Our construction leaves the function unaltered outside the compact set, preserves smoothness outside the singular set and convexity at a global scale. 
Furthermore, we prove that such singular perturbation can be locally and uniformly approximated with smooth convex functions in a neighborhood of the singular set. 
Once again, the difficulty here is to leave the functions unaltered outside the compact set.

\begin{lemma}
\label{lem:singular deformation}
Let $f: \R^n \to \R$ be a smooth, convex function. For every $x \in \R^n$, $r >0$ there exists a convex function $f_\infty: \R^n \to \R$ such that 
\begin{enumerate}[label=(\roman*)]
\item $f_\infty$ is smooth and equal to $f$  outside $B_r(x)$;
\item the graph of $f_\infty$ restricted to $B_r(x)$ has a dense set of singularities.
\end{enumerate}
Furthermore, there exists a sequence of smooth, strictly convex functions $f_\infty^k:\R^n \to \R$ converging uniformly to $f_\infty$ and equal to $f$ outside $B_r(x)$.
\begin{proof}
Take $\lbrace y_k\rbrace_{k = 1}^\infty$ any dense set contained in $S \coloneqq B_r(x)$.
We want to perturb $f$ in $S$ to obtain a new convex function whose graph has singularities in correspondence with $y_k$. 
To do so, we consider $g : B_1(0) \to \R$ such that
\begin{equation*}
    \begin{cases}
    g(x) = |x| + |x|^2 - 1  & x \in B_{1/2}(0)\\
    g \in C^\infty(B_1(0) \setminus \lbrace 0\rbrace) \\
    \supp g \subset B_1(0)\\
    g \leq 0.
    \end{cases}
\end{equation*}
Then, for $\varepsilon> 0$ and $y \in \R^n$ we define $g_{\varepsilon, y }: B_\varepsilon(y) \to \R$ as
\begin{equation*}
    g_{\varepsilon, y }(x) \coloneqq g\left(\frac{x-y}{\varepsilon} \right),
\end{equation*}
so that $g_{\varepsilon, y }$ is smooth outside $\lbrace y\rbrace$, non-positive and strictly convex on $B_{\varepsilon/2}(y)$.

Let $\varepsilon_1 < 1 - |y_1|$, define
\begin{equation*}
    f_1(x) \coloneqq f(x) + \eta_1 g_{\varepsilon_1, y_1}(x),
\end{equation*}
with $\eta_1>0$ small enough so that $f_1$ is strictly convex.
Observe that $f_1$ is smooth outside $\lbrace y_1 \rbrace$ and (its graph) has a singular point on $y_1$. 

Recursively, we let $\varepsilon_k < \min \left\lbrace 1-|y_k|, \text{dist} (y_k, y_1), \ldots \text{dist} (y_k, y_{k-1}) \right\rbrace$ and we define
\begin{equation}
    \label{eq: f_k}
    f_k(x) \coloneqq f_{k-1}(x) + \eta_k g_{\varepsilon_k, y_k}(x).
\end{equation}
By construction $f_k$ is smooth outside $\lbrace y_1, \ldots, y_k \rbrace$, where its graph is singular, and strictly convex, provided that $\eta_k$ is small enough.
Furthermore, if $\eta_k$ are such that $\sum_k \eta_k$ converges, $f_k$ converges uniformly to some $f_\infty$, which is convex, singular on $\lbrace y_k\rbrace_{k = 1}^\infty$ and is smooth elsewhere. Moreover it is equal to $f$ outside $S$.
Observe also that $\lbrace (y_k, f_\infty (y_k)) \rbrace_{k = 1}^\infty$ is dense in $\Graph(f_\infty|_S)$ since $f_\infty$ is locally Lipschitz.

It remains to show that $f_\infty$ can be smoothly approximated with strictly convex functions. 
By a diagonalization procedure, it is enough to uniformly approximate each $f_k$.

For $0 < \delta < \min\lbrace \varepsilon_1, \ldots, \varepsilon_k \rbrace$, let $\phi_{\delta, k} = \phi_\delta: \R^n \to \R$ be a smooth convex function such that
\begin{equation*}
    \phi_\delta(x) = f_k(x) \text{ on } \R^n \setminus \bigcup_{i=1}^k B_\delta(y_k).
\end{equation*}
The existence of $\phi_\delta$ is ensured by \cite[Theorem 2.1]{G2002}. 
Clearly $\phi_\delta$ converges pointwise to $f_k$ as $\delta \to 0$. Since the functions are all strictly convex the convergence is actually uniform. 
This concludes the proof. 
\end{proof}
\end{lemma}

\begin{remark}
Observe that the epigraph of $f_\infty$ is a convex set in $\R^{n+1}$ whose boundary, endowed with the intrinsic distance, is an Alexandrov space of non\hyp{}negative curvature (see \cite{B1976}). 
Its singularities are contained (and dense) in the compact set $\Graph(f_\infty|_S)$.
Similarly, the graphs of $f_\infty^k$ are smooth hypersurfaces of positive sectional curvature, isometrically immersed in $\R^{n+1}$. 
Since $f_\infty^k \to f_\infty$ uniformly, their graphs converge with respect to the Hausdorff metric.
In the case of convex sets of $\R^n$, it is well known that this implies Gromov\hyp{}Hausdorff convergence; see \cite[Theorem 10.2.6]{BBI2001} observing that the proof applies in any dimension.
Observe also that the convergence is measured if we endow these spaces with the usual $n$\hyp{}dimensional Hausdorff measure $\mathcal{H}^n$. 
On an isometrically immersed manifold, this is in fact the Riemannian volume. 
\end{remark}

\section{Convergence of solutions of the Poisson equation}\label{sec: Poisson}
The next step in our proof is a convergence result for the solutions of the Poisson equation on limit spaces. 
In what follows we mimic, up to minor modifications necessary to our purposes, \cite[Proposition B.1]{P2020}, where Pigola collects and develops a series of previous results due to Honda, \cite{H2015} and \cite{H2018}.

Let us consider the following space 
\begin{equation*}
    \mathcal{M}(n, D) = \lbrace (M, g) \text{ cpt.} : \dim M = n, \diam (M) \leq D, \Sect \geq 0 \rbrace, 
\end{equation*}
and denote with $\overline{\mathcal{M}(n, D)}$ its closure with respect to the measured Gromov\hyp{}Hausdorff topology. Note that elements of $\overline{\mathcal{M}(n, D)}$ are in particular Alexandrov spaces with $\Curv \geq 0$ and $\diam \leq D$. Note that, by volume comparison and bounds on the diameter, there exists $V > 0$ depending on $n, D$ such that $\vol X \leq V$ for all $X \in \mathcal{M}(n, D)$.

\begin{remark}
The following Proposition actually holds in the more general setting of Ricci limit spaces. To avoid unnecessary complication in notations, we restrict ourselves to the case of Alexandrov spaces which are a special case of the former. 
\end{remark}

In what follows, all convergences are intended in the sense of Honda, see \cite[Section 3]{H2015}.

\begin{proposition}
\label{prop:convergence poisson}
Let $(M_k, g_k) \in \mathcal{M}(n, D)$ be a sequence of smooth manifolds converging in the mGH topology to an Alexandrov space $(X_\infty, d_\infty,\mu_\infty) \in \overline{\mathcal{M}(n, D)}$ of dimension $n$ and let $x_\infty \in X_\infty$. There exist functions $u_k \in C^2(M_k)$, $g_k \in \Lip(M_k)$ and $u_\infty \in W^{1, 2}(X_\infty) \cap L^p(X_\infty)$, $g_\infty \in L^p(X_\infty)$ for all $1<p<+\infty$, such that $u_k, u_\infty$ are non-constant and $\Delta_{M_k} u_k = g_k$, $\Delta_{X_\infty} u_\infty = g_\infty$. Furthermore
\begin{enumerate}[label=\normalfont (\alph*)]
    \item $g_\infty \geq \sfrac{1}{2}$ on a neighborhood of $x_\infty$; 
    \item $g_k \to g_\infty$ in the strong $L^p$ sense;
    \item $u_k \to u_\infty$ in the strong $W^{1, 2}$ sense;
    \item $\Vert u_k \Vert_{W^{1, p}} \leq L$ for some $L = L(p, n, D, K) > 0$;
    \item $u_k \to u_\infty$ in the strong $L^p$ sense;
    \item $\nabla^{M_k} u_k \to \nabla^X u_\infty$ in the weak $L^p$ sense.
\end{enumerate}
These functions satisfy (a) through (f) for all $1<p<+\infty$.
\begin{proof}
Since $M_k$ is bounded, separable and $M_k$ converges to $X_\infty$ with respect to the mGH topology, there exists a sequence of points $x_k \in M_k$ such that the mGH convergence $(M_k, g_k, x_k) \to (X_\infty, \mu_\infty, x_\infty)$ is pointed.

Next, using volume comparison and the convergence $\vol (M_k)\to \mathcal H^n(X)$ as $k\to\infty$, one can show the existence of a uniform $R> 0$ such that for $k \gg 1$, 
\begin{equation*}
    \vol B_R^{M_k}(x_k) \leq \frac{1}{2} \vol M_k.
\end{equation*}
Let $f_k : M_k \to [0, 1]$ be Lipschitz functions compactly supported in $B_R^{M_k}(x_k)$ satisfying
\begin{equation*}
    i) \, f_k = 1\ \text{ on } B_{R/2}^{M_k}(x_k), \qquad ii) \, \Vert \nabla f_k \Vert_{L^\infty} \leq \frac{2}{R}.
\end{equation*}
Define 
\begin{equation*}
    g_k \coloneqq f_k - \fint_{M_k} f_k \in \Lip(M_k),
\end{equation*}
and note that
\begin{equation*}
    0 \leq \fint_{M_k} f_k \leq \frac{\vol B_R^{M_k}(x_k)}{\vol M_k} \leq \frac{1}{2}.
\end{equation*}
Clearly $\fint_{M_k} g_k = 0$ and $\Vert g_k\Vert_{L^\infty} \leq 1$. Moreover, $g_k \geq 1/2$ on $B_{R/2}^{M_k}(x_k)$ so that $g_k \not\equiv 0$. 
Since $\Vert g_k\Vert_{L^\infty} \leq 1$ and the volumes are uniformly bounded, $\Vert g_k \Vert_{L^p} \leq V^{1/p}$ for all $p > 1$.
Using \cite[Proposition 3.19]{H2015} we conclude that $g_k$ converges weakly to some $g_\infty \in L^p(X_\infty)$ (\cite[Definition 3.4]{H2015}).
Condition $ii)$ ensures that the sequence $g_k$ is asymptotically uniformly continuous in the sense of \cite[Definition 3.2]{H2015}. Hence, $g_k$ converges to $g_\infty$ strongly and in the sense of \cite[Definition 3.1]{H2015}, see \cite[Remark 3.8]{H2015}. 
This ensures that $g_\infty\not\equiv 0$ in a neighborhood of $x_\infty$ and, more importantly, allows us to use \cite[Proposition 3.32]{H2015} which proves strong $L^p$ convergence of $g_k$ to $g_\infty$. It is worthwhile to notice that $g_k$ converges $L^p$ strongly to $g_\infty$ for every $1<p<+\infty$, in particular, for $p = 2$. 

Next, we denote with $u_k \in C^2(M_k)$ the unique (non-constant) solution of the Poisson equation
\begin{equation*}
    \Delta_{M_k} u_k = g_k \quad \text{on }M_k, 
\end{equation*}
satisfying
\begin{equation*}
    \fint_{M_k} u_k = 0.
\end{equation*}
Since $g_k$ converges to $g_\infty$ in a strong (and thus weak) $L^2$ sense, \cite[Theorem 1.1]{H2018} ensures $W^{1, 2}$ convergence of $u_k$ to the unique (non-constant) solution $u_\infty \in W^{1, 2}(X_\infty)$ of the Poisson equation
\begin{equation*}
    \Delta_{X_\infty} u_\infty = g_\infty \quad \text{on }X_\infty, 
\end{equation*}
satisfying
\begin{equation*}
    \fint_{X_\infty} u_\infty = 0.
\end{equation*}

Finally, we claim that $\lbrace u_k \rbrace$ is bounded in $W^{1, p}$. By \cite[Theorem 4.9]{H2015} this implies $L^p$ strong convergence of $u_k$ to $u_\infty$ and $L^p$ weak convergence of $\nabla^{M_k} u_k$ to $\nabla^X u_\infty$ up to a subsequence and thus concludes the proof of \Cref{prop:convergence poisson}. 
To prove the claim we observe that since $u_k \to u_\infty$ in a strong $W^{1,2}$ sense, we have $L^2$ boundedness of $u_k$. Applying the estimates in \cite[Corollary 4.2]{ZZ2019} we obtain $L^\infty$ bounds for $ u_k$ and $\nabla u_k$, hence, the desired $L^p$ estimates using the uniform bound on volumes.
\end{proof}
\end{proposition}

\section{Proofs of the results}
\label{sec: proof}
In \Cref{sec:singular space} we established a method to locally perturb a smooth and strictly convex function by introducing a set of singular points, which is dense inside a given compact.
In the following we consider a sequence of infinitely many singular perturbations scattered over a non\hyp{}compact manifold, each of these perturbations is GH approximated with smooth Riemannian manifolds.  
For each perturbation, we prove that it is impossible to have the validity of a local (hence of a global) Calder\'on\hyp{}Zygmund inequality whose constant is uniformly bounded across the approximating sequence of manifolds. 
To do so, we show that each singular set together with its corresponding approximation can be seen as a piece of a compact space whose metric is smooth outside the singular part. 
This observation is a technical device which allows the application of already available results. In particular, we can employ \Cref{prop:convergence poisson} to localize the strategy of De Philippis and N\'u\~nez-Zimbr\'on in a neighborhood of each singular set.
Once we have proven that the constants of the local Calder\'on\hyp{}Zygmund inequalities cannot be chosen uniformly, we select on the $j$-th perturbation in the approximating sequence a manifold with $CZ(p)$ constant greater that $j$. 

\begin{lemma}
\label{lem:main lemma}
Let $n \geq 2$ and $p> n$. There exists a sequence of smooth and strictly convex functions $f_j : \R^n \to \R$, $j \geq 1$ and a monotone increasing sequence of radii $r_j > 0$ such that
\begin{enumerate}[label=(\roman*)]
\item $f_j(x) = f_{j-1}(x)$ for $x \in \mathfrak{B}_{j-1}$;
\item $f_j(x) = |x|^2$ for $x \in \R^n \setminus \overline{\mathfrak{B}_j}$;
\end{enumerate}
where $ \mathfrak{B}_j = B_{r_j}(0)$ and $\mathfrak{B}_0 = \emptyset$. 
Furthermore, if we consider $N_j = \Graph(f_j)$ as a Riemannian manifold isometrically immersed in $\R^{n+1}$, there exists some $v_j \in C^2(N_j)$ compactly supported in $\Graph (f_j|_{\mathfrak{B}_j \setminus \overline{\mathfrak{B}_{j-1}}}$) which satisfies
\begin{equation}
    \label{eq:contraddiction to CZ(p)}
    \Vert \Hess v_j \Vert_{L^p} > j \left(\Vert \Delta v_j \Vert_{L^p} + \Vert v_j \Vert_{L^p}\right),
\end{equation}
where $L^p = L^p(M_j)$.
\begin{proof}
We begin with a remark on notation: given a subset $A \subset \R^n$ and some function $k : \R^n \to \R$, we denote with $k(A) = \Graph(k|_A) \subset \R^{n+1}$.
This abuse of notation is repeatedly used throughout the proof.  

To simplify the exposition, the proof proceeds inductively on $j\ge 1$.
Set $f_0(x)=|x|^2$.
Suppose one has $f_{j-1}$ and wants to build $f_j$. 
Let $S_j$ be a Euclidean ball contained in $\R^n \setminus \overline{\mathfrak{B}_{j-1}}$. By \Cref{lem:singular deformation} there exists a convex function $h_j$ with a dense set of singular points in $S_j$ and equal to $f_{j-1}$ outside $S_j$.  
Furthermore, $h_j$ can be approximated with smooth and strictly convex functions $h_{j,k} : \R^n \to \R$ equal to $f_{j-1}$ outside $S_j$.
Note that $h_j(S_j)$ corresponds to the $\mathfrak{D_j}$ of the Introduction.

Next, let $r_j > 0$ be such that $S_j \subset \mathfrak{B}_j$. For later use we observe that one can always consider a larger ball $T_j$ such that $S_j \subset T_j$ and $T_j \subset \mathfrak{B}_j \setminus \overline{\mathfrak{B}_{j-1}}$.
We want to extend $h_j(\mathfrak{B}_j)$ to a closed (i.e. compact without boundary) Alexandrov space $X_j$ with $\Curv(X_j) \geq 0$. Moreover, we would like the extension to be smooth outside $h_j(\mathfrak{B}_{j-1})$. 
To this purpose, let $A_j$ be the upper hemisphere of boundary $h_j(\partial \mathfrak{B}_j)$ in $\R^{n+1}$, so that $\widetilde{X_j} \coloneqq    h_j(\mathfrak{B}_j) \cup A_j$ is a convex hypersurface in $\R^{n+1}$. 
To obtain $X_j$, one simply needs to smooth $\widetilde{X_j}$ in a neighborhood of $h_j(\partial \mathfrak{B}_j)$. 
For instance, one can use \cite[Theorem 2.1]{G2002}, observing that in this neighborhood, $\widetilde{X_j}$ is obtained by rotation of a piecewise smooth curve.
We consider on $X_j$ the metric induced by $\R^{n+1}$. 
By the same strategy, we extend $h_{j,k}(\mathfrak{B}_j)$ to a compact and smooth Riemannian manifold $M_{j,k}$ with $\Sect(M_{j,k}) > 0$, isometrically immersed in $\R^{n+1}$.

Note that, for all $k$, $ M_{j,k}=X_j$ outside of $S_j$.  Moreover $M_{j,k}$ converges to $X_j$ in a (measured) Gromov\hyp{}Hausdorff sense as $k \to \infty$.
Then, choosing a point $x_{j,\infty} \in  S_j \subset X_j$, we apply \Cref{prop:convergence poisson} to deduce the existence of $u_{j,k} \in C^2(M_{j,k})$, $g_{j,k} \in \Lip(M_{j,k})$ and $u_{j,\infty} \in W^{1, 2}(X_j) \cap L^p(X_j)$, $g_{j,\infty} \in L^p(X_j)$ such $\Delta_{M_{j,k}} u_{j,k} = g_{j,k}$ and $\Delta_{X_j} u_{j,\infty} = g_{j,\infty}$.
In particular
\begin{enumerate}[label=(\alph*)]
    \item $\Delta u_{j,k} \to \Delta u_{j,\infty}$ strongly in $L^p$, hence, $\Vert \Delta u_{j,k}\Vert_{L^p} \leq C_1$;
    \item $\Vert u_{j,k}\Vert_{W^{1,p}} \leq C_1$;
    \item $g_{j,\infty} \geq 1/2$ in a neighborhood of $x_{j,\infty}$. 
    In particular, in this neighborhood $u_{j,\infty}$ can not be constant.
\end{enumerate}
Here $C_1$ depends on $n, p$ and the upper bound $\diam M_{j,k} \leq D_j$ and the norms are intended over $L^p=L^p(M_{j,k})$ and $W^{1,p} = W^{1,p}(M_{j,k})$.

A key element in our proof is the possibility to localize the sequence $u_{j,k}$ without altering its essential properties. 
This can be done via smooth cutoff functions $\chi_{j,k} \in C^\infty(M_{j,k})$ equal to $1$ on $h_{j,k}(S_j)$ and identically $0$ outside of $h_{j,k}(T_j)$. 
Moreover, since the manifolds $M_{j,k}$ are all isometric outside $h_{j,k}(S_j)$, we can choose the functions $\chi_j = \chi_{j,k}$ so that they are equal (independently of $k$) outside $h_{j,k}(S_j)$.
Let $v_{j,k} \coloneqq \chi_j\, u_{j,k} \in C^2(M_{j,k})$ and observe that $v_{j,k}$ preserves the $L^p$ bounds of $u_{j,k}$, indeed: 
\begin{equation}
\label{eq:v_k}
    \Vert v_{j,k} \Vert_{L^p} \leq \Vert u_{j,k} \Vert_{L^p}\leq C_2,
\end{equation}
\begin{equation} 
\label{eq:delta v_k}
    \Vert \Delta v_{j,k} \Vert_{L^p} \leq \Vert \Delta u_{j,k} \Vert_{L^p} + \Vert u_{j,k} \Delta \chi_j \Vert_{L^p} + 2 \Vert |\nabla u_{j,k}|\,|\nabla \chi_j| \Vert_{L^p} \leq C_2, 
\end{equation}
where $C_2$ depends on $C_1$ as well as on the choice of $\chi_j$.

Next, we need some function theoretic considerations.  
First, we observe that compactness of $M_{j,k}$ implies the validity of an $L^p$\hyp{}Calder\'on\hyp{}Zygmund inequality
\begin{equation}
    \label{eq:CZ(p)}
    \Vert \Hess \varphi \Vert_{L^p} \leq E_{j,k} \left(\Vert \Delta \varphi \Vert_{L^p} + \Vert \varphi\Vert_{L^p}\right), \quad \forall\,\varphi \in C^2(M_{j,k}).
\end{equation}
Second, if $p > n$, we have the validity on the sequence $M_{j,k}$ of a uniform Morrey\hyp{}Sobolev inequality 
\begin{equation}
\label{eq:morrey embedding}
     \left\vert \varphi(x) - \varphi(y) \right\vert \leq C_3 \Vert \nabla \varphi \Vert_{L^p} d_{j,k}(x, y)^{1 - \frac{n}{p}}, \quad \forall\,\varphi \in C^1(M_{j,k}), 
\end{equation}
where $d_{j,k}$ is the Riemannian distance on $M_{j,k}$, and the constant $C_3$ depends on $n, p$ and the uniform upper bound on $\diam M_{j,k}$.
See \cite[Theorem 9.2.14]{HK2015} for reference, observing that the lower bound on the Ricci curvature ensures the validity of a $p$-Poincaré inequality; see \cite[Theorem 5.6.5]{S2002}.
Applying \eqref{eq:morrey embedding} to $|\nabla \varphi|$ and combining with the Calder\'on\hyp{}Zygmund inequality \eqref{eq:CZ(p)} implies the following estimate
\begin{equation}
    \label{eq: morrey + CZ(p)}
    \left\vert |\nabla \varphi |(x) - |\nabla \varphi|(y) \right\vert \leq C_3 E_{j,k} \left(\Vert \Delta \varphi \Vert_{L^p} + \Vert \varphi \Vert_{L^p}\right) d_{j,k}(x, y)^{1 -\frac{n}{p}}, 
\end{equation}
for all $\varphi \in C^2(M_{j,k})$ and all $x, y \in M_{j,k}$.

Applying \eqref{eq: morrey + CZ(p)} to $v_{j,k}$ and using estimates \eqref{eq:v_k} and \eqref{eq:delta v_k} we obtain
\begin{equation}
\label{eq:asympt unif cont v_k}
     \left\vert |\nabla v_{j,k} |(x) - |\nabla v_{j,k}|(y) \right\vert \leq C E_{j,k} d_{j,k}(x, y)^{1-\frac{n}{p}} \quad x, y \in M_{j,k}, 
\end{equation}
where $C$ depends on $C_1, C_2$ and $C_3$, i.e., $C = C(n, p, \chi_j, D_j)$.
Suppose by contradiction that $E_{j,k}$ is bounded from above uniformly in $k$.
By \eqref{eq:asympt unif cont v_k} we deduce that $|\nabla v_{j,k}|$ is uniformly asymptotic continuous in the sense of Honda, hence, from \cite[Proposition 3.3]{H2015} we conclude that $|\nabla v_{j,k}|$ converges pointwise to $|\nabla v_{j,\infty}| \in C^0(X)$. 
However, since $X$ is an $n$\hyp{}dimensional Alexandrov space with $\Sect \geq 0$, it is a $RCD(0, n)$ space. Moreover, $\Delta v_{j,\infty} \in L^{p>n}$. 
From \cite[Theorem 1.1]{DPNZ2019} we then conclude that $|\nabla v_{j,\infty} |(x) = |\nabla u_{j,\infty}|(x) = 0$ whenever $x$ is a singular point. 
Note here that singular points of Alexandrov spaces are \textit{sharp} in the sense of De Philippis and N\'u\~nez-Zimbr\'on and have finite Bishop-Gromov density.
By density we conclude that, $v_{j,\infty}$ must be constant in a neighborhood of $x_{j,\infty}$ thus contradicting (c). 

In particular there exists some $\bar{k}$, which may depend on $j$, such that
\begin{equation}
    \label{eq:CZ(p) v_k}
     \Vert \Hess v_{j,\bar{k}} \Vert_{L^p} > j \left(\Vert \Delta v_{j,\bar{k}} \Vert_{L^p} + \Vert v_{j,\bar{k}}\Vert_{L^p}\right)
\end{equation}
on $M_{j,\bar{k}}$.
Finally, we set $f_j = h_{j,\bar{k}}$, since $v_{j,\bar{k}}$ is compactly supported in $h_{j,\bar{k}}(T_j)$, it defines a function $v_j = v_{j,\bar{k}}$ on $N_j$ which satisfies \eqref{eq:contraddiction to CZ(p)}. 
\end{proof}
\end{lemma}

Note that, while \Cref{prop:convergence poisson} is independent of $p$, the previous result depends on the initial choice of $p>n$. 
This has to be attributed to the fact that the constants $C_1, C_2, C_3$ and $C$ are all dependent on $p$.

To obtain a contradiction to $CZ(p)$ for $p>n$, we then simply need to glue the manifolds of \Cref{lem:main lemma} together.

\begin{myproof}[of \Cref{t: main}]
For $p>n$, let $f_j$ be as in \Cref{lem:main lemma}, and let $f$ be its point-wise limit. Note that the convergence is actually uniform on compact sets. 
The function $f$ is smooth and strictly convex, thus, $M = \Graph(f)$ is a smooth, non\hyp{}compact Riemannian manifold isometrically immersed in $\R^{n+1}$ satisfying $\Sect(M) > 0$. 
Since $f$ is defined on the whole space $\R^n$, $M$ is also a complete manifold. 
Observe that the sequence $v_j$ as in \Cref{lem:main lemma} induces functions in $C^2(M)$ whose supports are compact and disjoint, and which satisfy \eqref{eq:contraddiction to CZ(p)} on $L^p(M)$. 
This sequence clearly contradicts the validity of a global Calder\'on\hyp{}Zygmund inequality on $M$.
\end{myproof}

Note that in the above we have not exploited to the fullest the fact that the functions $v_j$ have disjoint supports. 
In fact, not only one has a sequence $v_j$ on which \eqref{eq:contraddiction to CZ(p)} holds, but one can actually define a function $F \in C^2(M)$ such that $\Vert F \Vert_{L^p} + \Vert \Delta F \Vert_{L^p} < + \infty$ but $\Vert \Hess F \Vert_{L^p} = + \infty$, which is a stronger condition. 
This allows to prove \Cref{cor:counterexample to sobolev inclusion}.

\begin{myproof}[of \Cref{cor:counterexample to sobolev inclusion}]
Fix $p >n $, let $(M,g)$ and $v_j \in C^2(M)$ be as in the proof of \Cref{t: main}.
Define
\begin{equation*}
    F \coloneqq \sum_{j=1}^\infty \frac{1}{j^2}\frac{v_j}{\Vert \Delta v_j \Vert_{L^p} + \Vert v_j \Vert_{L^p}},
\end{equation*}
and observe that the sum converges since it is locally finite. 
Note that
\begin{equation*}
    \Vert \Delta F \Vert_{L^p} + \Vert F \Vert_{L^p}=\sum_{j=1}^\infty \frac{1}{j^2},
\end{equation*}
so that $F \in H^{2, p}(M)$.
By \eqref{eq:contraddiction to CZ(p)}, on the other hand, we have
\begin{equation*}
    \Vert \Hess F \Vert_{L^p} \geq \sum_{j=1}^\infty \frac{1}{j}, 
\end{equation*}
hence, $F \not\in W^{2, p}(M)$.
\end{myproof}

We conclude this paper with a proof of Theorem \ref{th: sharp}
which follows quite directly from the constructions of the counterexamples in \cite{GP2015,L2020}. 
These counterexamples rely on the construction of manifolds whose sectional curvature are increasingly oscillating on a sequence of compact annuli going to infinity. 
However, by distancing the (disjoint) annuli far enough we are able to provide a controlled lower bound on sectional curvatures.   
\begin{myproof}[of \Cref{th: sharp}]
The counterexamples to \eqref{e: CZp} in \cite{GP2015,L2020} are constructed on a model manifold $(M, g)$, i.e. $M=[0,+\infty)\times \mathbb S^{n-1}$ endowed with a warped metric $g = dt^2 +\sigma^2(t) g_{\mathbb S^{n-1}}$.
By carefully choosing the warping function $\sigma$, the authors proved the existence of a sequence of smooth functions $\{u_k\}_{k=1}^\infty$ and a sequence of intervals $\{[a_k,b_k]\}_{k=1}^\infty$ such that 
    \begin{itemize}
        \item $a_{k+1}>b_k$;
        \item $u_k$ is compactly supported in the annulus $[a_k,b_k] \times \mathbb S^{n-1}$;
        \item the sequence of functions $u_k$ contradicts \eqref{e: CZp} for any possible constant, i.e.  \[
        \frac{\Vert \Hess u_k \Vert_{L^p}}{\Vert \Delta u_k \Vert_{L^p} + \Vert u_k \Vert_{L^p}}\to \infty,\qquad\text{as }k\to\infty;
        \]
        \item  there exists two sequences of intervals $\{[c_k,d_k]\}_{k=1}^\infty$ and $\{[e_k,f_k]\}_{k=1}^\infty$ with $b_k<c_k<d_k<e_k<f_k<a_{k+1}$ such that $\sigma$ is linear and increasing on $[c_k,d_k]$ and is linear and decreasing on $[e_k,f_k]$, namely
        \[\sigma|_{[c_k,d_k]}(t)= \alpha_k t + \beta_k,\quad\text{and}\quad \sigma|_{[e_k,f_k]}(t)= \gamma_k t + \delta_k
        \]
        for some constants $\alpha_k>0$, $\gamma_k <0$ and $\beta_k,\delta_k\in\R$.
        \end{itemize}
Note that, in order to satisfy this latter condition, our $\{u_k\}_{k=1}^\infty$ could be a subsequence of the sequence $\{u_k\}_{k=1}^\infty$ produced in \cite{L2020}

Now, for $k\ge 2$, let $0<\kappa_k<\infty$ be such that 
\[
\forall\, x\in [e_{k-1},d_k]\times \mathbb S^{n-1},\quad \min\Sect(x)\ge 
-\kappa_k.
\]
Up to an increase of $\kappa_{k+1}$, we can assume that $\kappa_k\le \kappa_{k+1}$.
For $k\ge 2$, let $T_k$ be such that $\lambda (T_k)> \kappa_k$.
For later purpose, since $\lambda$ is increasing we can assume without loss of generality that $T_{k+1}>T_{k}+d_{k-1}-e_{k-2}$ and that 
\begin{equation}
\label{condition T}
    \alpha_{k-1}(T_{k+1}+e_{k-2}-T_k)+\beta_{k-1}>\sigma(e_{k-1}).
\end{equation}
We define now a new warping function $\tilde\sigma(t):[0,+\infty)\to [0,+\infty)$ and a corresponding model metric $\tilde g=dt^2 +\tilde \sigma^2(t) g_{\mathbb S^{n-1}}$ on $M$ as follows.
We define $\tilde \sigma(t)$ only for $t\ge T_3$, since the choice of $\tilde \sigma$ on $[0,T_3)$ does not affect the conclusion of the theorem. For $t\in [T_k,T_k+d_{k-1}-e_{k-2}]$ define 
\[\tilde\sigma(t)=\sigma(t+e_{k-2}-T_k),\]
so that 
\[
\Sect_{\tilde g}\ge -\kappa_{k-1} 
\]
on $[T_k,T_k+d_{k-1}-e_{k-2}]\times \mathbb S^{n-1}$. In particular
\[
\Sect_{\tilde g}(t,\Theta)\ge -\kappa_{k}> - \lambda (T_k)\ge -\lambda (t) 
\]
for any $(t,\Theta)\in ([T_k,T_k+d_{k-1}-e_{k-2}]\cup [T_{k+1},T_{k+1}+d_{k}-e_{k-1}])\times \mathbb S^{n-1}$. It remains to prescribe $\tilde \sigma$ on the intervals $(T_k+d_{k-1}-e_{k-2},T_{k+1})$ for $k\ge 3$.
Note that on $[T_k+c_{k-1}-e_{k-2},T_k+d_{k-1}-e_{k-2}]$ we have $\tilde\sigma(t)=\alpha_{k-1}(t+e_{k-2}-T_k)+\beta_{k-1}$. Similarly, on $[T_{k+1},T_{k+1}+f_{k-1}-e_{k-1}]$, we have 
$\tilde\sigma(t)=\gamma_{k-1}(t+e_{k-1}-T_{k+1})+\delta_{k-1}$. Because of assumption \eqref{condition T}, we can find a $S_k\in (T_k+d_{k-1}-e_{k-2},T_{k+1})$ such that 
\[\hat \sigma (t)=\begin{cases}\alpha_{k-1}(t+e_{k-2}-T_k)+\beta_{k-1}&\text{on }[T_k+c_{k-1}-e_{k-2},S_k]\\\gamma_{k-1}(t+e_{k-1}-T_{k+1})+\delta_{k-1}&\text{on }
[S_k,T_{k+1}+f_{k-1}-e_{k-1}]
\end{cases}\]
is a well-defined concave continuous piece-wise linear function which coincides with $\tilde\sigma$ outside  $(T_k+d_{k-1}-e_{k-2},T_{k+1})$. Let $\epsilon_k>0$ be a small constant to be fixed later, and define $\tilde\sigma$ on  $(T_k+d_{k-1}-e_{k-2},T_{k+1})$ to be a concave smooth approximation of $\hat\sigma$ equal to $\hat\sigma$ outside $[S_k-\epsilon_k,S_k+\epsilon_k]$ (this can be produced for instance applying \cite[Theorem 2.1]{G2002}).
A standard computation show that the sectional curvature of $(M,\tilde g)$ are given by
\[
\Sect_{rad} (t,\Theta) = -\frac{\tilde \sigma '' (t)}{\tilde\sigma (t)},\qquad \Sect_{tg} (t,\Theta) = \frac{1-(\tilde \sigma ' (t))^2}{\tilde\sigma (t)^2},\]
for tangent planes respectively containing the radial direction, or orthogonal to it. 
Since $\tilde \sigma$ is concave for  $t \in (T_k+d_{k-1}-e_{k-2},T_{k+1})$ then
\[\Sect_{rad} (t,\Theta)\ge 0 \ge -\lambda (t).\]
If $\alpha_{k-1}\leq 1$ and  $\gamma_{k-1}\geq -1$ then $\Sect_{tg}(t,\Theta) \geq 0 \geq -\lambda(t)$ in a trivial way. 
Otherwise, 
\[
\Sect_{tg}(t,\Theta)>\Sect_{tg}(T_k+d_{k-1}-e_{k-2},\Theta)\ge -\kappa_k>-\lambda (t)
\] 
for $t \in (T_k+d_{k-1}-e_{k-2},S_k-\epsilon_k)$ and 
\[\Sect_{tg}(t,\Theta)>\Sect_{tg}(T_{k+1},\Theta)\ge -\kappa_k>-\lambda (t)\] 
for $t \in (S_k+\epsilon_k, T_{k+1})$.
Finally, for $t\in[S_k-\epsilon_k,S_k+\epsilon_k]$, by concavity 
\[
1-(\tilde\sigma'(t))^2\ge \min\{1-(\tilde\sigma'(S_k-\epsilon_k))^2;1-(\tilde\sigma'(S_k+\epsilon_k))^2\},\]
while $\tilde\sigma(t)$ is arbitrarily close to $\tilde\sigma(S_k-\epsilon_k)$ and to $\tilde\sigma(S_k+\epsilon_k)$ for $\epsilon_k$ small enough. 
Accordingly, we can choose $\epsilon_k$ small enough so that $\Sect_{tg}(t,\Theta)>-\lambda(t)$ also for $t \in [S_k-\epsilon_k,S_k+\epsilon_k]$. 
Hence, we have proved that for all $t\ge T_3$, the sectional curvature of $(M,\tilde g)$ at $(t,\Theta)$ are lower bounded by $-\lambda(t)$. 
Observe that
$([T_k,T_k+d_{k-1}-e_{k,2}]\times\mathbb S^{n-1},\tilde g)$ is isometric to $([e_{k,2},d_{k-1}]\times\mathbb S^{n-1},g)$.
Then we conclude by defining  $w_k(t,\Theta)=u_{k-1}(t+e_{k-2}-T_{k},\Theta)$ so that the $w_k$ are smooth, compactly supported in $[T_k+a_{k-1}-e_{k-2},T_k+b_{k-1}-e_{k-2}]\times\mathbb S^{n-1}$ and verify
\[
    \frac{\Vert \Hess w_k \Vert_{L^p}}{\Vert \Delta w_k \Vert_{L^p} + \Vert w_k \Vert_{L^p}}\to \infty,\qquad\text{as }k\to\infty.
    \]
\end{myproof}

\bibliography{main}
\bibliographystyle{acm}
\end{document}